\RequirePackage{fix-cm}
\documentclass[smallextended]{svjour3}       
\smartqed  
\usepackage{latexsym,amssymb,amsmath}
\usepackage{amssymb}
\usepackage{amsmath}
\usepackage{graphicx}
 \journalname{Bull. Math. Sci.}
\begin{document}

\title{Geometry of warped product semi-slant submanifolds of Kenmotsu manifolds
}
%



\author{Siraj Uddin
}


\institute{S. Uddin \at
              Department of Mathematics,\\
               Faculty of Science,\\
                King Abdulaziz University, 21589 Jeddah, Saudi Arabia\\
              \email{siraj.ch@gmail.com}           
}

\date{Received: date / Accepted: date}
\sloppy
\maketitle

\begin{abstract}
In this paper, we study semi-slant submanifolds and their warped products in Kenmotsu manifolds. The existence of such warped products in Kenmotsu manifolds is shown by an example and a characterization. A sharp relation is obtained as a lower bound of the squared norm of second fundamental form in terms of the warping function and the slant angle. The equality case is also considered in this paper. Finally, we provide some applications of our derived results. 
\end{abstract}

\section{Introduction}
\label{intro}

Warped product manifolds were defined and studied by Bishop and OÕNeill in 1969 as a natural generalization of the Riemannian product manifolds \cite{Bi}. The geometrical aspects of these manifolds have been studied later by many mathematicians.

The study of warped product submanifolds from extrinsic view points was initiated by B.-Y. Chen in \cite{C2,C3}. He proved that there do not exist warped product CR-submanifolds of the form $M_\perp \times_f M_T$, where $M_T$ and $M_\perp$ are holomorphic (complex) and totally real submanifolds of a Kaehler manifold $\tilde M$, respectively. Then he studied warped products of the form $M_T \times_f M_\perp$, known as CR-warped products. Several fundamental results on CR-warped products in Kaehler manifolds were established in  \cite{C2,C3,Book} by him. Motivated by ChenÕs fundamental seminal work, many geometers studied warped product submanifolds in various Riemannian manifolds (see \cite {Al}, \cite{C4}, \cite{Ha}, \cite{Mi}, \cite{Mun}, \cite{U2} and the references therein). For the most recent detailed survey on warped product manifolds and submanifolds, see \cite{book17}.

The notion of slant submanifolds of almost Hermitian manifolds was introduced by Chen in \cite{C1}.  Later, J.L. Cabrerizo et al. studied  in \cite{Ca2} slant immersions in K-contact and Sasakian manifolds. In particular, they provided interesting examples of slant submanifolds in both almost contact metric manifolds and Sasakian manifold. In \cite{Ca2}, they also characterized slant submanifolds by means of the covariant derivative of the square of the tangent projection on the submanifold. In  \cite{Ca1}, they defined  and studied semi-slant submanifolds of Sasakian manifolds.

The non-existence of warped product semi-slant submanifolds of Kaehler manifolds was proved in \cite{S1}. Later, M. Atceken studied warped product semi-slant submanifolds and proved non-existence of such submanifolds in Kenmotsu manifolds \cite{At1}. The warped product semi-slant submanifolds of Kenmotsu manifolds were also studied in \cite{Mus,U1}. For the survey on warped product submanifolds of Kenmotsu manifolds, we refer to \cite{At2,Mur,Mus,Ol}.

In this paper, we continue the study of warped product semi-slant submanifolds in Kenmotsu manifolds. In the first part of this paper, we give some fundamental results for semi-slant submanifolds of Kenmotsu manifolds. Then we prove existence of such warped products by applying a characterization result. We also give an example of non-trivial proper warped product semi-slant submanifolds in Kenmotsu manifolds. In the second part, we derive a sharp inequality for the squared norm of the second fundamental form in terms of the warping function and the slant angle. We also investigate the equality case of this inequality. Several applications of our results are present in the last part.

\section{Preliminaries}

Let $\tilde M$ be an almost contact metric manifold with
structure $(\varphi, \xi, \eta, g)$, where $\varphi$ is a $(1,1)$ tensor
field, $\xi$ a vector field, $\eta$ is a $1$-form and $g$ is a
Riemannian metric on $ \tilde M$ satisfying the following properties
\begin{align}
\label{2.1}
\varphi^2=-I+\eta\otimes\xi,~~\varphi\xi=0,~~\eta\circ\varphi=0,~~\eta(\xi) = 1.
\end{align}
\begin{align}
\label{2.2}
g(\varphi X, \varphi Y)=g(X, Y)-\eta(X)\eta(Y).
\end{align}
In addition, if the following relation
\begin{align}
\label{2.3}
(\tilde\nabla_X\varphi)Y=g(\varphi X, Y)\xi-\eta(Y)\varphi X
\end{align}
holds for any $X, Y$ on $\tilde M$, then $\tilde M$ is called a \textit{Kenmotsu manifold}, where $\tilde\nabla$ is the Levi-Civita connection of $g$. It is easy to see from \eqref{2.3} that $\tilde\nabla_X\xi=X-\eta(X)\xi$. We shall use the symbol $\Gamma(T\tilde M)$ for the Lie algebra of vector fields on the manifold $\tilde M$.

Let $M$ be a submanifold of an almost contact metric manifold $\tilde M$ with induced metric $g$ and if $\nabla$ and $\nabla^{\perp}$ are the induced connections on the tangent bundle $TM$ and the normal bundle $T^{\perp}M$ of $M$, respectively then Gauss-Weingarten formulas are respectively given by
\begin{align}
\label{2.4}
\tilde\nabla_XY=\nabla_XY+h(X, Y)
\end{align}
\begin{align}
\label{2.5}
\tilde\nabla_XN=-A_NX+\nabla^{\perp}_XN,
\end{align}
for each $X,~Y\in\Gamma(TM)$ and $N\in\Gamma(T^\perp M)$, where $h$ and $A_N$ are the second
fundamental form and the shape operator (corresponding to the normal vector
field $N$) respectively for the immersion of $M$ into $\tilde M$. They are related by the relation $g(A_NX, Y)=g(h(X, Y), N)$.

For any $X\in\Gamma(TM)$ and $N\in\Gamma(T^\perp M)$ , we write
\begin{align}
\label{2.6}
(a)\,\,\varphi X=PX+FX,\,\,\,\,(b)\,\,\varphi N=tN+fN
\end{align}
where $PX$ and $tN$ are the tangential components of $\varphi X$ and $\varphi N$, respectively and $FX$ and $fN$ are the normal components of $\varphi X$ and $\varphi N$, respectively.

A submanifold $M$ is said to be {\it invariant} if $F$ is identically zero, that is, $\varphi X\in\Gamma(TM)$ for any $X\in\Gamma(TM)$. On the other hand, $M$ is {\it anti-invariant} if $P$ is identically zero, that is, $\varphi X\in\Gamma(T^\perp M)$, for any $X\in\Gamma(TM)$.

Let $M$ be a submanifold tangent to the structure vector field $\xi$ isometrically immersed into an almost contact metric manifold $\tilde M$. Then $M$ is said to be a contact CR-submanifold if there exists a pair of orthogonal distributions ${\mathcal{D}}:p\to {\mathcal {D}}_p$ and ${\mathcal{D}}^\perp:p\to {\mathcal{D}}^\perp_p$,~$\forall~p\in M$ such that
\begin{enumerate}
\item [(i)] $TM=\mathcal D\oplus\mathcal D^\perp\oplus\langle\xi\rangle$, where $\langle\xi\rangle$ is the $1$-dimensional distribution spanned by the structure vector field $\xi$,
\item [(ii)] $\mathcal D$ is invariant, i.e., $\varphi\mathcal D=\mathcal D$,
\item [(iii)] $\mathcal D^\perp$ is anti-invariant, i.e., $\varphi{\mathcal D^\perp}\subseteq T^\perp M$.
\end {enumerate}

The invariant and anti-invariant submanifolds are the special cases of a contact CR-submanifold. If we denote the dimensions of $\mathcal D$ and ${\mathcal D^\perp}$ by $d_1$ and $d_2$, respectively then $M$ is {\it{invariant}} (resp. {\it{anti-invariant}}) if $d_2=0$ (resp. $d_1=0$). 

There is another class of submanifolds which are known as slant submanifolds which we define as follows:

For each non zero vector $X$ tangent to $M$ at $p$, such that
$X$ is not proportional to $\xi$, we denote by $\theta(X)$, the angle
between $\varphi X$ and $T_pM,$ for all $p\in M$. Then, $M$ is said to be \textit{slant} \cite{Ca2} if the angle $\theta(X)$
is constant for all $X\in TM-\{\xi\}$ and $p\in M$ i.e., $\theta(X)$ is independent of the choice of the vector field $X$ and the point $p\in M$. The angle
$\theta(X)$ is called the \textit{slant angle}.
Obviously, if $\theta=0$ then, $M$ is invariant and if $\theta=\pi/2$ then, $M$
is an anti-invariant submanifold. If the slant angle of $M$ is
neither 0 nor $\pi/2$, then it is called \textit{proper
slant}.

A characterization of slant submanifolds was given in \cite{Ca2} as follows:

\begin{theorem} \cite{Ca2} Let $M$ be a submanifold of an
almost contact metric manifold $\tilde M$ such that $\xi\in\Gamma(TM)$. Then $M$ is
slant if and only if there exists a constant $\lambda\in[0,1]$ such that
\begin{align}
\label{2.7}
P^2=\lambda(-I+\eta\otimes\xi).
\end{align}
Furthermore, in such case, if $\theta$ is slant angle, then $\lambda=\cos^2\theta$.
\end{theorem}

 The following relations are straight forward consequences of \eqref{2.7}
\begin{equation}
\label{2.8}
g(PX,PY)=\cos ^{2}\theta \lbrack g(X,Y)-\eta (X)\eta (Y)]\
\end{equation}
\begin{equation}
\label{2.9}
g(FX,FY)=\sin ^{2}\theta \lbrack g(X,Y)-\eta (X)\eta (Y)]
\end{equation}
for any $X,Y$ tangent to $M.$

\section{Semi-slant submanifolds}

In \cite{P}, semi-slant submanifolds were defined and studied by N. Papaghiuc as a natural generalization of CR-submanifolds of almost Hermitian manifolds in terms of the slant distribution.  Later on, Cabrerizo et al. \cite{Ca1} studied these submanifolds in contact geometry. They defined these submanifolds as follows:

\begin{definition} \cite{Ca1} \rm{A Riemannian submanifold $M$ of an almost contact manifold $\tilde M$ is said to be a {\it{semi-slant submanifold}} if there exist two orthogonal distributions ${\mathcal{D}}$ and ${\mathcal{D}}^\theta$ such that $TM={\mathcal{D}}\oplus {\mathcal{D}}^\theta\oplus\langle\xi\rangle$,~the distribution ${\mathcal{D}}$ is invariant  i.e., $\varphi{\mathcal{D}}={\mathcal{D}}$ and the distribution ${\mathcal{D}}^\theta$ is slant with slant angle $\theta\neq \frac{\pi}{2}.$}
\end{definition}

If we denote the dimensions of ${\mathcal D}$ and ${\mathcal{D}}^\theta$  by $d_1$ and $d_2$ respectively, then it is clear that contact CR-submanifolds and slant submanifolds are semi-slant submanifolds with $\theta=\frac{\pi}{2}$ and $d_1=0$, respectively. If neither $d_1=0$ nor $\theta=\frac{\pi}{2}$, then $M$ is a \textit{proper semi-slant} submanifold.

Moreover, if $\nu$ is the $\varphi-$invariant subspace of the normal bundle $T^\perp M$, then in case of a semi-slant submanifold, the normal bundle $T^\perp M$ can be decomposed as $T^\perp M= F{\mathcal D}^\theta\oplus \nu.$

First, we give the following non-trivial example of a semi-slant submanifold of an almost contact metric manifold.

\begin{example} \rm{Consider a submanifold $M$ of ${\mathbb{R}^9}$ with the cartesian coordinates $(x_1,\,x_2,\,x_3,\,x_4,\,y_1,\,y_2,\,y_3,\, y_4,\,z)$ and the contact structure
\begin{align*}
\varphi\left(\frac{\partial}{\partial x_i}\right)=-\frac{\partial}{\partial y_i},~~~~\varphi\left(\frac{\partial}{\partial y_j}\right)=\frac{\partial}{\partial x_j},~~~~\varphi\left(\frac{\partial}{\partial z}\right)=0,~~~~1\leq i, j\leq4.
\end{align*}
It is easy to show that $(\varphi, \xi, \eta, g)$ is an almost contact metric structure on ${\mathbb{R}^9}$ with $\xi=\frac{\partial}{\partial z},\,\eta=dz$ and $g$, the Euclidean metric of ${\mathbb{R}^9}$. Consider an immersion $\psi$ on ${\mathbb{ R}}^9$ defined by
\begin{align*}
&\psi(\theta,\,\phi,\,v,\,w,\,z)\\
&=(\cos(\theta+\phi),\,\theta-\phi,\,\frac{1}{2}\theta+\phi,\,v+w,\,\sin(\theta+\phi),\,\phi-\theta,\,\theta+\frac{1}{2}\phi,\,w-v,\,z).
\end{align*}
If we put 
\begin{align*}
X_1=-\sin(\theta+\phi)\frac{\partial}{\partial x_1}+\frac{\partial}{\partial x_2}+\frac{1}{2}\frac{\partial}{\partial x_3}+\cos(\theta+\phi)\frac{\partial}{\partial y_1}-\frac{\partial}{\partial y_2}+\frac{\partial}{\partial y_3},
\end{align*}
\begin{align*}
X_2=-\sin(\theta+\phi)\frac{\partial}{\partial x_1}-\frac{\partial}{\partial x_2}+\frac{\partial}{\partial x_3}+\cos(\theta+\phi)\frac{\partial}{\partial y_1}+\frac{\partial}{\partial y_2}+\frac{1}{2}\frac{\partial}{\partial y_3},
\end{align*}
\begin{align*}
X_3=\frac{\partial}{\partial x_4}-\frac{\partial}{\partial y_4},\,\,X_4=\frac{\partial}{\partial x_4}+\frac{\partial}{\partial y_4},\,\,X_5=\frac{\partial}{\partial z}
\end{align*}
then the restriction of $\{X_1,X_2,X_3,X_4, X_5\}$ to $M$ forms an orthogonal frame fields of the tangent bundle $TM$. Clearly, we have
\begin{align*}
\varphi X_1=\sin(\theta+\phi)\frac{\partial}{\partial y_1}-\frac{\partial}{\partial y_2}-\frac{1}{2}\frac{\partial}{\partial y_3}+\cos(\theta+\phi)\frac{\partial}{\partial x_1}-\frac{\partial}{\partial x_2}+\frac{\partial}{\partial x_3},
\end{align*}
\begin{align*}
\varphi X_2=\sin(\theta+\phi)\frac{\partial}{\partial y_1}+\frac{\partial}{\partial y_2}-\frac{\partial}{\partial y_3}+\cos(\theta+\phi)\frac{\partial}{\partial x_1}+\frac{\partial}{\partial x_2}+\frac{1}{2}\frac{\partial}{\partial x_3},
\end{align*}
\begin{align*}
\varphi X_3=-\frac{\partial}{\partial y_4}-\frac{\partial}{\partial x_4},\,\,\varphi X_4=-\frac{\partial}{\partial y_4}+\frac{\partial}{\partial x_4},\,\,\varphi X_5=0.
\end{align*}
It is easy to verify that $\mathcal D=\rm{Span}\{X_3,\,X_4\}$ is an invariant distribution and $\mathcal D^\theta=\rm{Span}\{X_1,\,X_2\}$ is a slant distribution of $M$ with slant angle $\theta=\cos^{-1}\left(\frac{3}{17}\right)$ such that $X_5=\xi=\frac{\partial}{\partial z}$ is tangent to $M$. Thus, $M$ is a proper semi-slant submanifold of ${\mathbb{R}^9}$.}
\end{example}

Now, we give some basic results for semi-slant submanifolds of Kenmotsu manifolds which are useful to the next section.

\begin{lemma} Let $M$ be a proper semi-slant submanifold of a Kenmotsu manifold $\tilde M$ with invariant and proper slant distributions ${\mathcal{D}}\oplus<\xi>$ and ${\mathcal D}^\theta$, respectively. Then we have:
\begin{align}
\label{3.1}
\sin^2\theta\,g(\nabla_YX, Z)=g(A_{FZ}\varphi X-A_{FPZ}X, Y)
\end{align}
for any $X, Y\in\Gamma({\mathcal{D}}\oplus\langle\xi\rangle)$ and $Z\in\Gamma({\mathcal{D}}^\theta)$.
\end{lemma}
\begin{proof} For any $X, Y\in\Gamma({\mathcal{D}}\oplus\langle\xi\rangle)$ and $Z\in\Gamma({\mathcal{D}}^\theta)$, we have
\begin{align*}
g(\nabla_YX, Z)=g(\tilde\nabla_YX, Z)=g(\varphi\tilde\nabla_YX, \varphi Z).
\end{align*}
Using the covariant derivative property of $\varphi$ and the relation \eqref{2.6}, we get
\begin{align*}
g(\nabla_YX, Z)=g(\tilde\nabla_Y\varphi X, PZ)+g(\tilde\nabla_Y\varphi X, FZ)-g((\tilde\nabla_Y\varphi)X, \varphi Z).
\end{align*}
Then from \eqref{2.3}, \eqref{2.4} and the fact that $\varphi X$ and $PZ$ are orthogonal vector fields, we derive
\begin{align*}
g(\nabla_YX, Z)&=-g(\varphi X, \tilde\nabla_YPZ)+g(h(\varphi X, Y), FZ)\notag\\
&=g(X, \tilde\nabla_Y\varphi PZ)-g(X, (\tilde\nabla_Y\varphi)PZ)+g(A_{FZ} \varphi X, Y)\notag\\
&=g(X, \tilde\nabla_YP^2Z)+g(X, \tilde\nabla_YFPZ)-g((\tilde\nabla_Y\varphi)PZ, X)+g(A_{FZ} \varphi X, Y).\notag
\end{align*}
Again using \eqref{2.3},  \eqref{2.5} and  \eqref{2.7} and the fact that $\xi\in\Gamma({\mathcal{D}})$, we find
\begin{align*}
g(\nabla_YX, Z)=-\cos^2\theta g(X, \tilde\nabla_YZ)-g(A_{FPZ}X, Y)+g(A_{FZ} \varphi X, Y).
\end{align*}
Hence, the result follows from the last relation. 
\end{proof}

The following corollary is an immediate consequence of the above lemma.
\begin{corollary} Let $M$ be a proper semi-slant submanifold of a Kenmotsu manifold  $\tilde M$. Then, the distribution ${\mathcal{D}}\oplus\langle\xi\rangle$ defines a totally geodesic foliation if and only if
\begin{align*}
g(A_{FZ}\varphi X-A_{FPZ}X, Y)=0
\end{align*}
for any $X, Y\in\Gamma({\mathcal{D}}\oplus\langle\xi\rangle)$ and $Z\in\Gamma({\mathcal{D}}^\theta)$.
\end{corollary}

Also, we have the following results for the leaves of the slant distribution.
\begin{lemma} Let $M$ be a proper semi-slant submanifold of a Kenmotsu manifold $\tilde M$. Then we have
\begin{align}
\label{3.2}
g(\nabla_ZW, X)=\csc^2\theta\,\big(g(A_{FPW}X-A_{FW}\varphi X, Z)\big)-\eta(X)g(Z, W)
\end{align}
for any $X\in\Gamma({\mathcal{D}}\oplus\langle\xi\rangle)$ and $Z, W\in\Gamma({\mathcal{D}}^\theta)$.
\end{lemma} 
\begin{proof} For any $X\in\Gamma({\mathcal{D}}\oplus\langle\xi\rangle)$ and $Z, W\in\Gamma({\mathcal{D}}^\theta)$, we find
\begin{align*}
g(\nabla_ZW, X)=g(\tilde\nabla_Z\varphi W,\varphi X)+\eta(X)g(\tilde\nabla_ZW, \xi).
\end{align*}
Using \eqref{2.3} and \eqref{2.6}, we get
\begin{align*}
g(\nabla_ZW, X)=g(\tilde\nabla_ZPW,\varphi X)+g(\tilde\nabla_ZFW,\varphi X)-\eta(X)g(Z, W).
\end{align*}
Then from \eqref{2.2} and \eqref{2.5}, we arrive at
\begin{align*}
g(\nabla_ZW, X)=-g(\varphi\tilde\nabla_ZPW, X)-g(A_{FW}Z,\varphi X)-\eta(X)g(Z, W).
\end{align*}
By using the covariant derivative property of $\varphi$ and the symmetry of the shape operator $A$, we obtain
\begin{align*}
g(\nabla_ZW, X)=g((\tilde\nabla_Z\varphi)PW, X)-g(\tilde\nabla_Z\varphi PW, X)-g(A_{FW}\varphi X, Z)-\eta(X)g(Z, W).
\end{align*}
Again using \eqref{2.3} and \eqref{2.6}, we derive
\begin{align*}
g(\nabla_ZW, X)&=-g(\tilde\nabla_ZP^2W, X)-g(\tilde\nabla_ZFPW, X)+g(PZ, PW)\eta(X)\notag\\
&-g(A_{FW}\varphi X, Z)-\eta(X)g(Z, W).
\end{align*}
From Theorem 1 and the relation \eqref{2.5}, we find
\begin{align*}
g(\nabla_ZW, X)&=\cos^2\theta g(\tilde\nabla_ZW, X)+g(A_{FPW}Z, X)+\cos^2\theta g(Z, W)\eta(X)\notag\\
&-g(A_{FW}\varphi X, Z)-\eta(X)g(Z, W).
\end{align*}
Since $A$ is self-adjoint, then by using trigonometric identities, the result follows from the last relation.
\end{proof}

From the above lemma, we have the following results.
\begin{corollary} Let $M$ be a proper semi-slant submanifold of a Kenmotsu manifold  $\tilde M$. Then the slant distribution ${\mathcal{D}}^\theta$ defines a totally geodesic foliation if and only if
\begin{align*}
g(A_{FPZ}X-A_{FZ}\varphi X, W)=\eta(X)g(Z,W)
\end{align*}
for any $X\in\Gamma({\mathcal{D}}\oplus\langle\xi\rangle)$ and $Z, W\in\Gamma({\mathcal{D}}^\theta)$.
\end{corollary}
\begin{lemma} Let $M$ be a proper semi-slant submanifold of a Kenmotsu manifold $\tilde M$. Then we have:
\begin{align*}
\sin^2\theta\,g([Z, W], X)=g(A_{FZ}\varphi X-A_{FPZ}X, W)-g(A_{FW}\varphi X-A_{FPW}X, Z)
\end{align*}
for any $X\in\Gamma({\mathcal{D}}\oplus\langle\xi\rangle)$ and $Z, W\in\Gamma({\mathcal{D}}^\theta)$.
\end{lemma}
\begin{proof} From \eqref{3.2}, we have
\begin{align}
\label{3.3}
\sin^2\theta\,g(\tilde\nabla_ZW, X)=g(A_{FPW}X-A_{FW}\varphi X, Z)-\sin^2\theta\,\eta(X)g(Z, W)
\end{align}
for any $X\in\Gamma({\mathcal{D}}\oplus\langle\xi\rangle)$ and $Z, W\in\Gamma({\mathcal{D}}^\theta)$. Interchanging $Z$ and $W$ in \eqref{3.3}, we find
\begin{align}
\label{3.4}
\sin^2\theta\,g(\tilde\nabla_WZ, X)=g(A_{FPZ}X-A_{FZ}\varphi X, W)-\sin^2\theta\,\eta(X)g(Z, W).
\end{align}
Thus, the result follows from \eqref{3.3} and \eqref{3.4}.
\end{proof}

\section{Warped product semi-slant submanifolds}
In this section, we study warped product semi-slant submanifolds of Kenmotsu manifolds, by considering that one of the factor is a slant submanifold. In the following first, we give brief introduction of warped product manifolds.

Let $M_{1}$ and $M_{2}$ be two Riemannian manifolds with Riemannian metrics $g_{1}$ and $g_{2},$ respectively, and $f$ be a positive differentiable function on $M_{1}$. Then,  $M=M_{1}\times _{f}M_{2}=(M_{1}\times M_{2}, g),$ is a warped product manifold of $\ M_{1}$ \ and $M_{2}$ such that
\begin{align}
\label{4.1}
g(X, Y)=g_1({\pi_1}_\star X, {\pi_1}_\star Y)+(f\circ\pi_1)^2g_2({\pi_2}_\star X, {\pi_2}_\star Y)
\end{align}
where the vector fields $X$ and $Y$ are tangent to $
M=M_{1}\times _{f}M_{2}$ at $(p,q)$ and $\pi _{1}$ and $\pi _{2}$ are the canonical projections of $M=M_{1}\times M_{2}$ onto $M_{1}$ and $M_{2}$, respectively and $\star$ is the symbol for the tangent map. A warped product manifold $M=M_{1}\times _{f}M_{2}=(M_{1}\times M_{2},g)$ is said to be {\it{trivial}} or simply a \textit{{Riemannian product}} if the warping function $f$ is constant. If $X$ is a vector field on $M_1$ and $V$ is an another vector field on $M_2$, then from Lemma 7.3 of \cite{Bi}, we have
\begin{equation}
\label{main}
\nabla_XV=\nabla_VX=X(\ln f)V
\end{equation}
where $\nabla$ denotes the Levi-Civita connection on $M$.  It is well-known that $M_1$  is a totally geodesic submanifold and $M_2$ is a totally umbilical submanifold of  $M$ (cf. \cite{Bi,Book}). 
 
 The gradient $\vec\nabla f$ of a function $f$ on $M$ is defined as $g(\vec\nabla f,X)=X(f),$ for any vector field $X$ on $M$. If $\{e_1,\cdots e_n\}$ is an orthonormal frame field of the tangent space of $M$, then we have
\begin{align}
\label{4.2}
{\Vert \vec\nabla f \Vert}^{2}=\sum_{i=1}^{n}{(e_{i}(f))}^{2}.
\end{align}

In this paper, we study warped product semi-slant submanifolds of the form $M_T\times_fM_\theta$ of a Kenmotsu manifold $\tilde M$, where $M_T$ and $M_\theta$ are invariant and proper slant submanifolds of $\tilde M$, respectively. First, we consider  $M_1$ and $M_2$ be two Riemannian submanifolds of a Kenmotsu manifold $\tilde M$. Then their warped product submanifold $M$ is of the form $M_1\times_fM_2$. Since, we consider the structure vector field $\xi$ is tangent to $M$, therefore two possible cases arise:\\

{\textit{Case 1}}. When the structure vector field $\xi$ is tangent to $M_2$, then for any $X\in\Gamma(TM_1)$, we have $\tilde\nabla_X\xi=X$. From \eqref{2.4} and \eqref{main}, we find that $X(\ln f)\xi=X$, by taking the inner product with $\xi$, we observe that $f$ is constant, i.e., the warped product $M_1\times_fM_2$ becomes a Riemannian product (trivial).\\

{\textit{Case 2}}. When the structure vector field $\xi$ is tangent to $M_1$, then for any $Z\in\Gamma(TM_2)$, we have
\begin{align*}\tilde\nabla_Z\xi=\nabla_Z\xi+h(Z, \xi).
\end{align*}
Then from \eqref{2.3} and \eqref{main}, we find that
\begin{align}
\label{4.3}
(i)~~\xi\ln f=1,~~~~~~~~~~~~~~~~~~~~~~~~(ii)~~h(Z, \xi)=0,~~\forall~Z\,\in\Gamma(TM_2).\end{align}

Now, in the following we consider the warped products of the form $M_T\times_fM_\theta$, called warped product semi-slant submanifolds of a Kenmotsu manifold $\tilde M$ such that $\xi$ is tangent to $M_T$, where $M_T$ and $M_\theta$ are invariant and proper slant submanifolds of $\tilde M$, respectively. If neither $\dim M_T$ is zero nor the slant angle of $M_\theta$ is $\frac{\pi}{2}$, then the warped product semi-slant submanifold is called proper. It is clear that the contact CR-warped product submanifolds are the special cases of warped product semi-slant submanifolds.

First, we give the following non-trivial example of warped product semi-slant submanifolds in Kenmotsu manifolds.

\begin{example} \rm{Consider the complex Euclidean space $\mathbb C^4$ with its usual Kaehler structure and the real global coordinates $(x_1, x_2, x_3, x_4, y_1, y_2, y_3, y_4)$. Let $\tilde M=\mathbb R\times_f\mathbb C^4$ be a warped product manifold between the product of real line $\mathbb R$ and the complex space $\mathbb C^4$. In fact, $\tilde M$ is a Kenmotsu manifold with the almost contact metric structure $(\varphi, \xi, \eta, g)$ such that
\begin{align*}\varphi\left(\frac{\partial}{\partial x_i}\right)=-\frac{\partial}{\partial y_i},~~~\varphi\left(\frac{\partial}{\partial y_j}\right)=\frac{\partial}{\partial x_j},~~~\varphi\left(\frac{\partial}{\partial z}\right)=0,~~~1\leq i, j\leq 4,
\end{align*}
and
\begin{align*}\xi=e^z\left(\frac{\partial}{\partial z}\right),~~~\eta=e^zdz,~~~~~g=e^{2z}<, >\end{align*}
where $<, >$ denotes the Euclidean metric tensor of $\mathbb R^9$. Consider a submanifold $M$ defined by the immersion $\phi$ as follows
\begin{align*}
\phi(u_1,u_2,u_3,u_4,z)=(u_1,0,u_3,0,u_2, 0, u_4\cos\theta, u_4\sin\theta, z)
\end{align*}
with $\theta\in\left(0, \pi/2\right)$. Then the tangent space $TM$ of $M$ at any point is spanned by the following vectors
\begin{align*}
X_1=\frac{\partial}{\partial x_1}+y_1\left(\frac{\partial}{\partial z}\right),~~X_2=\frac{\partial}{\partial y_1},~~X_3=\frac{\partial}{\partial x_3}+y_3\left(\frac{\partial}{\partial z}\right),\end{align*}
\begin{align*}
X_4=\cos\theta\frac{\partial}{\partial y_3}+\sin\theta\frac{\partial}{\partial y_4},~~X_5=\frac{\partial}{\partial z}=\frac{1}{e^z}\xi.
\end{align*}
Then, we find
\begin{align*}
\varphi X_1=-\frac{\partial}{\partial y_1},~~\varphi X_2=\frac{\partial}{\partial x_1},~~\varphi X_3=-\frac{\partial}{\partial y_3},
\end{align*}
\begin{align*}\varphi X_4=\cos\theta\frac{\partial}{\partial x_3}+\sin\theta\frac{\partial}{\partial x_4},~~\varphi X_5=0.\end{align*}
Thus, $M$ is a proper semi-slant submanifold tangent to the structure vector field $\xi$ with invariant and proper slant distributions $\mathcal{D}=\rm{Span}\{X_1, X_2\}$ and $\mathcal{D}^\theta=\rm{Span}\{X_3, X_4\}$ respectively with slant angle $\theta$. It is easy to see that the distributions ${\mathcal{D}}$ and $\mathcal{D}^\theta$ are integrable. Consider, the integral manifolds corresponding to the distributions $\mathcal{D}$ and $\mathcal{D}^\theta$ by $M_T$ and $M_\theta$, respectively. Then it is easy to check that $M=M_T\times_fM_\theta$ is a proper warped product semi-slant submanifold isometrically immersed in $\tilde M$ with warping function $f=e^z,\,z\in\mathbb R$.}
\end{example} 

Now, we have the following useful lemma for later use.

\begin{lemma} Let $M=M_{T}\times {_{f}}M_{\theta }$
be a warped product semi-slant submanifold of a Kenmotsu manifold
$\tilde{M}$ such that $\xi$ is tangent to $M_{T}$, where $M_{T}$ is an invariant submanifold and
$M_{\theta }$ is a proper slant submanifold of
$\tilde{M}$. Then, we have
\begin{align}
\label{4.4}
g(h(X, Z), FW)=\big(\eta(X)-X(\ln f)\big)g(Z, PW)-\varphi X(\ln f)\,g(Z, W)
\end{align}
for any $X\in\Gamma(TM_T)$ and $Z, W\in\Gamma(TM_\theta)$.
\end{lemma}

\begin{proof} For any $X\in\Gamma(TM_T)$ and $Z, W\in\Gamma(TM_\theta)$, we have

\begin{align*}
g(h(X, Z), FW)&=g(\tilde\nabla_ZX, \varphi W)-g(\tilde\nabla_ZX, PW)\notag\\
&=-g(\varphi\tilde\nabla_ZX, W)-g(\nabla_ZX, PW).
\end{align*}
Using a covariant derivative property of $\varphi$ and \eqref{main}, we find
\begin{align*}
g(h(X, Z), FW)=g((\tilde\nabla_Z\varphi)X, W)-g(\tilde\nabla_Z\varphi X, W)-X(\ln f)\,g(Z, PW).
\end{align*}
Again using \eqref{2.3}, \eqref{2.4} and \eqref{main}, we get the desired result.
\end{proof}

The following relations can be easily obtained by interchanging $X$ by $\varphi X$ and $Z$ by $PZ$ and $W$ by $PW$ in \eqref{4.4}, for any $X\in\Gamma(TM_T)$ and $Z, W\in\Gamma(TM_\theta)$
\begin{align*}
g(h(\varphi X, Z), FW)=\big(X(\ln f)-\eta(X)\big)g(Z, W)-\varphi X(\ln f)\,g(Z, PW),\end{align*}
\begin{align}
\label{4.5}
g(h(X, PZ), FW)&=\varphi X(\ln f)\,g(Z, PW)\notag\\
&-\cos^2\theta\big(X(\ln f)-\eta(X)\big)g(Z, W),
\end{align}
\begin{align}
\label{4.6}
g(h(X, Z), FPW)&=\cos^2\theta\big(X(\ln f)-\eta(X)\big)g(Z, W)\notag\\
&-\varphi X(\ln f)\,g(Z, PW),\end{align}
and
\begin{align}
\label{4.7}
g(h(X, PZ), FPW)&=-\cos^2\theta\varphi X(\ln f)g(Z, W)\notag\\
&-cos^2\theta\big(X(\ln f)-\eta(X)\big)\,g(Z, PW).
\end{align}
From \eqref{4.5} and \eqref{4.6}, we have
\begin{align*}
g(h(X, PZ), FW)=-g(h(X, Z), FPW).
\end{align*}

On the other hand, for any $X, Y\in\Gamma(TM_T)$ and $Z\in\Gamma(TM_\theta)$, we have by Lemma 4.1 (i) of \cite{Mus}
\begin{align}
\label{4.8}
g(h(X, Y), FZ)=0.
\end{align}

In order to give a characterization result for semi-slant submanifolds of a Kenmotsu manifold, we recall the following result of Hiepko \cite{Hi}:\\

\noindent{\bf{Hiepko's Theorem.}} {\rm{Let ${\mathcal{D}}_1$ and ${\mathcal{D}}_2$ be two orthogonal distribution on a Riemannian manifold $M$. Suppose that ${\mathcal{D}}_1$ and ${\mathcal{D}}_2$ both are involutive such that ${\mathcal{D}}_1$ is a totally geodesic foliation and ${\mathcal{D}}_2$ is a spherical foliation. Then $M$ is locally isometric to a non-trivial warped product $M_1\times_fM_2$, where $M_1$ and $M_2$ are integral manifolds of ${\mathcal{D}}_1$ and ${\mathcal{D}}_2$, respectively.}}

Now, we are able to prove the following main result of this section.

\begin{theorem} Let $M$ be a proper semi-slant submanifold with invariant distribution $\mathcal{D}\oplus<\xi >$ and proper slant distribution $\mathcal{D}^{\theta }$ of a Kenmotsu manifold $\tilde{M}$. Then $M$ is locally a warped
product submanifold of the form $M_T\times_fM_\theta$ if and only if
\begin{align}
\label{4.9}
A_{FZ}\varphi X-A_{FPZ}X=\sin^2\theta\big(X(\mu)-\eta (X)\big)Z
\end{align}
for any $X\in\Gamma(\mathcal{D}\oplus<\xi>)$ and $Z\in\Gamma(\mathcal{D}^{\theta})$ and for some
function $\mu$ on $M$ satisfying $W\mu=0$, for any $W\in\Gamma(\mathcal{D}^{\theta })$.
\end{theorem}
\begin{proof} If $M=M_{T}\times _{f}M_{\theta}$ is a warped product submanifold of a Kenmotsu manifold $\tilde{M}$ such that $M_T$ is an invariant submanifold and $M_\theta$ is a proper slant submanifold of $\tilde M$, then from \eqref{4.8}, we have $g(A_{FZ}{\varphi X}, Y)=0$, for any $X, Y\in\Gamma(TM_T)$ and $Z\in\Gamma(TM_\theta)$, i.e., $A_{FZ}{\varphi X}$ has no components in $TM_T$. Also, if we interchange $Z$ by $PZ$ in \eqref{4.8}, then we have $g(A_{FPZ}{X}, Y)=0$, i.e., $A_{FPZ}{X}$ also has no components $TM_T$. Therefore, $A_{FZ}{\varphi X}-A_{FPZ}{X}$ lies in $TM_\theta$ only. On the other hand, for any $X, Y\in\Gamma(TM_T)$ and $Z, W\in\Gamma(TM_\theta)$, we have
\begin{align*}
g(A_{FZ}\varphi X, W)=g(h(\varphi X, W), FZ)=g(\tilde\nabla_W\varphi X, FZ).
\end{align*}
From the covariant derivative proper of $\varphi$, we find
\begin{align*}
g(A_{FZ}\varphi X, W)=g((\tilde\nabla_W\varphi)X, FZ)+g(\varphi\tilde\nabla_WX, \varphi Z)-g(\varphi\tilde\nabla_WX, PZ).
\end{align*}
Using \eqref{2.2}, \eqref{2.3}, \eqref{2.6}, \eqref{2.9} and \eqref{main}, we obtain
\begin{align*}
g(A_{FZ}\varphi X, W)&=-\sin^2\theta\,\eta(X)\,g(Z, W)+X(\ln f)g(Z, W)+g(\tilde\nabla_WX, P^2Z)\notag\\
&+g(\tilde\nabla_WX, FPZ).
\end{align*}
Then from \eqref{2.4}, Theorem 1 and \eqref{main}, we find
\begin{align*}
g(A_{FZ}\varphi X, W)&=-\sin^2\theta\,\eta(X)\,g(Z, W)+X(\ln f)g(Z, W)\notag\\
&-\cos^2\theta\,X(\ln f)\,g(Z, W)+g(A_{FPZ}X, W).
\end{align*}
Then \eqref{4.9} follows from the above relation with $\mu=\ln f$.

Conversely, if $M$ is a proper semi-slant submanifold of a Kenmotsu manifold $\tilde{M}$ such that \eqref{4.9} holds, then from Lemma 1 and the given condition \eqref{4.9}, we conclude that $\sin^2\theta\,g(\nabla_YX, Z)=0$, for any $X, Y\in\Gamma({\mathcal{D}}\oplus<\xi>)$ and $Z\in\Gamma({\mathcal{D}}^\theta)$. Since $M$ is a proper semi-slant submanifold, then $\sin^2\theta\neq0$, therefore $g(\nabla_YX, Z)=0$, i.e., the leaves of the distribution $\mathcal{D}\oplus <\xi>$ are totally geodesic in $M$. Also, from Lemma 3 and the given condition \eqref{4.9}, we find that $\sin^2\theta\,g([Z, W], X)=0$, for any $X\in\Gamma({\mathcal{D}}\oplus<\xi>)$ and $Z, W\in\Gamma({\mathcal{D}}^\theta)$. Since, $M$ is a proper semi-slant submanifold, thus we have $g([Z, W], X)=0$, i.e., the slant distribution ${\mathcal{D}}^\theta$ is integrable. If we consider $h^\theta$ be the second fundamental form of a leaf $M_\theta$ of ${\mathcal{D}}^\theta$ in $M$, then for any $Z, W\in\Gamma({\mathcal{D}}^\theta)$ and $X\in\Gamma({\mathcal{D}}\oplus<\xi>)$, we have
\begin{align*}
g(h^\theta(Z, W), X)=g(\nabla_ZW, X).\end{align*}
Using Lemma 2, we derive
\begin{align*}
g(h^\theta(Z, W), X)=-\csc^2\theta\,g(A_{FW}\varphi X-A_{FPW}X, Z)-\eta(X)g(Z, W).
\end{align*}
From the given condition \eqref{4.9}, we find
\begin{align*}
g(h^\theta(Z, W), X)=-X(\mu)g(Z,W).
\end{align*}
Then, from the definition of gradient, we get
\begin{align*}
h^\theta(Z, W)=-\vec\nabla(\mu)\,g(Z,W),
\end{align*}
which means that $M_\theta$ is totally umbilical in $M$ with mean curvature vector $H^\theta=-\vec\nabla(\mu)$. Now we prove that $H^\theta$ is parallel corresponding to the normal
connection $D^N$ of $M_{\theta }$ in $M$. Consider any $Y\in\Gamma(\mathcal{D}\oplus <\xi>)$ and $Z\in\Gamma(\mathcal{D}^{\theta})$, then we find that $g(D^N _{Z}\vec\nabla\mu ,Y)=g(\nabla _{Z}\vec\nabla\mu ,Y)=Zg(\vec\nabla\mu ,Y)-g(\vec\nabla\mu ,\nabla _{Z}Y)=Z(Y(\mu))-g(\vec\nabla\mu , [Z, Y])-g(\vec\nabla\mu ,\nabla _{Y}Z)=Y(Z\mu)+g(\nabla _{Y}\vec\nabla\mu ,Z)\}=0,$ since $Z(\mu )=0,\,\text{for all}\,Z\in\Gamma(\mathcal{D}^{\theta})$ and thus $\nabla _{Y}\vec\nabla\mu\in\Gamma(\mathcal{D}\oplus <\xi >)$. This means that the mean
curvature of $M_{\theta }$ is parallel. Thus the leaves of the distribution $\mathcal{D}^{\theta }$ are totally umbilical in $M$ with non-vanishing parallel mean curvature vector $H^\theta$ i.e., $M_\theta$ is an extrinsic sphere in $M$. Then from Hiepko's Theorem \cite{Hi}, $M$ is a warped product submanifold, which proves the theorem completely.
\end{proof}

Now, we construct the following orthonormal frame fields for a proper warped product semi-slant submanifold $M=M_T\times_fM_\theta$ of a $(2m+1)$-dimensional Kenmotsu manifold $\tilde M$ such that the structure vector field $\xi$ is tangent to $M_T$. Consider $M=M_T\times_fM_\theta$ be an $n$-dimensional warped product semi-slant submanifold of a Kenmotsu manifold $\tilde M$. If $\dim M_T=2t+1$ and $\dim M_\theta=2s$, then $%
n=2t+1+2s$. Let us consider the orthonormal frame fields of the corresponding
tangent bundles $\mathcal{D}$ and $\mathcal{D}^\theta$ of $M_T$ and $M_\theta$, respectively as: $\{e_1,\cdots, e_t, e_{t+1}=\varphi e_1,\cdots, e_{2t}=\varphi
e_t, e_{2t+1}=\xi\}$ is the orthonormal frame field of $\mathcal{D}$ and $%
\{e_{2t+2}=e_1^{\star},\cdots, e_{2t+1+s}=e_{s}^{\star}, e_{2t+2+s}=e_{s+1}^{\star}=\sec\theta\,Pe_{1}^{\star},\cdots, e_{n}=e_{2s}^{\star}=\sec\theta\,Pe_{s}^{\star}\}$
is the orthonormal frame field of $\mathcal{D}^\theta$. Then the orthonormal frame fields in the
normal bundle $T^\perp M$ of $F\mathcal{D}^\theta$ and invariant normal subbundle $\nu$ respectively are  $\{e_{n+1}=\tilde e_1=\csc\theta\,Fe_1^\star\cdots, e_{n+s}=\tilde
e_s=\csc\theta\,Fe_s^\star, e_{n+s+1}=\tilde e_{s+1}=\csc\theta\sec\theta\,FPe_1^\star,\cdots, e_{n+2s}=\tilde e_{2s}=\csc\theta\sec\theta\,FPe_s^\star\}$
and $\{e_{n+2s+1}=\tilde e_{2s+1},\cdots, e_{2m+1}=\tilde e_{2m+1-n-2s}\}$.

Next, we use the above constructed frame fields to find a relation (lower bound) for the squared norm of the second fundamental form $\|h\|^2$, in terms of the warping function and the slant angle of a proper warped product semi-slant submanifold of Kenmotsu manifolds.

\begin{theorem} Let $M=M_T\times{_f}M_\theta$ be a proper warped product semi-slant submanifold of a Kenmotsu manifold $\tilde M$ such that the structure vector field $\xi$ is
tangent to $M_T$, where $M_T$ is an invariant submanifold and $M_\theta$ is a proper slant submanifold of $\tilde M$. Then
\begin{enumerate}
\item[(i)] ~\textit{The squared norm of the second fundamental form $h$ of $M$ satisfies}
\begin{align}
\label{4.10}
\|h\|^2\geq4s\left(\csc^2\theta+\cot^2\theta\right)\left(\|\nabla^T(\ln f)\|^2-1\right)
\end{align}
where $\nabla^T(\ln f)$ is the gradient of the warping function $\ln f$ along $M_T$ and $2s=\dim M_\theta$.
\item[(ii)] If the equality sign in (i) holds, then $M_T$ is a
totally geodesic submanifold and $M_\theta$ is a totally umbilical submanifold of $\tilde M$. Furthermore, $M$ is minimal in $\tilde M$.
\end{enumerate}
\end{theorem}
\begin{proof}  From the definition of $h$, we have
\begin{align}
\label{4.11}
\|h\|^2&=\sum_{i, j=1}^n g(h(e_i, e_j), h(e_i, e_j))\notag\\
&=\sum_{r=n+1}^{2m+1}\sum_{i, j=1}^{n} g(h(e_i, e_j), e_r)^2\notag\\
&=\sum_{r=1}^{2s}\sum_{i, j=1}^{n} g(h(e_i, e_j), \tilde e_r)^2+\sum_{r=2s+1}^{2m+1-n-2s}\sum_{i, j=1}^{n}g(h(e_i, e_j), \tilde e_r)^2.
\end{align}
First term in the right hand side of \eqref{4.11} has $F\mathcal{D}^\theta$-components and the second term has $\nu$-components. Let us compute $F\mathcal{D}^\theta$-components terms only by using the frame fields of $\mathcal{D}$ and $\mathcal{D}^\theta$. Then we have
\begin{align}
\label{4.12}
\|h\|^2&\geq\sum_{r=1}^{2s} \sum_{i, j=1}^{2t+1} g(h(e_i, e_j), \tilde e_r)^2+2\sum_{r=1}^{2s}\sum_{i=1}^{2t+1}\sum_{ j=1}^{2s} g(h(e_i, e_j^\star), \tilde e_r)^2\notag\\
&+\sum_{r=1}^{2s} \sum_{i, j=1}^{2s} g(h(e_i^\star, e_j^\star), \tilde e_r)^2.
\end{align}
Then from \eqref{4.8}, the first term in the right hand side of above inequality
is identically zero. Let us compute just next term
\begin{align*}
\|h\|^2&\geq2\sum_{r=1}^{2s}\sum_{i=1}^{2t+1}\sum_{j=1}^{2s}g(h(e_i, e_j^\star), \tilde e_r)^2.\notag\\
&=2\sum_{r, j=1}^{2s}\sum_{i=1}^{2t}g(h(e_i, e_j^\star), \tilde e_r)^2+2\sum_{r, j=1}^{2s}g(h(\xi, e_j^\star), \tilde e_r)^2.
\end{align*}
Since for a submanifold $M$ of a Kenmotsu manifold $\tilde M$, we have $h(X, \xi)=0$, for any $%
X\in\Gamma(TM)$. By using this fact the last term in above inequality is
identically zero. Then from the assumed frame fields of $\mathcal{D}$, $\mathcal{D}^\theta$ and $%
F\mathcal{D}^\theta$, we derive
\begin{align*}
\|h\|^2&\geq2\csc^2\theta\sum_{i=1}^{t}\sum_{r, j=1}^{s}g(h(e_i, e_j^\star), Fe_r^\star)^2\notag\\
&+2\csc^2\theta\sec^4\theta\sum_{i=1}^{t}\sum_{r, j=1}^{s}g(h(e_i, Pe_j^\star), FPe_r^\star)^2\notag\\
&+2\csc^2\theta\sum_{i=1}^{t}\sum_{r, j=1}^{s}g(h(\varphi e_i, e_j^\star), Fe_r^\star)^2\notag\\
&+2\csc^2\theta\sec^4\theta\sum_{i=1}^{t}\sum_{r, j=1}^{s}g(h(\varphi e_i, Pe_j^\star), FPe_r^\star)^2\notag\\
&+2\csc^2\theta\sec^2\theta\sum_{i=1}^{t}\sum_{r, j=1}^{s}g(h(e_i, Pe_j^\star), Fe_r^\star)^2\notag\\
&+2\csc^2\theta\sec^2\theta\sum_{i=1}^{t}\sum_{r, j=1}^{s}g(h(e_i, e_j^\star), FPe_r^\star)^2\notag\\
&+2\csc^2\theta\sec^2\theta\sum_{i=1}^{t}\sum_{r, j=1}^{s}g(h(\varphi e_i, Pe_j^\star), Fe_r^\star)^2\notag\\
&+2\csc^2\theta\sec^2\theta\sum_{i=1}^{t}\sum_{r, j=1}^{s}g(h(\varphi e_i, e_j^\star), FPe_r^\star)^2.
\end{align*}
From the relations \eqref{4.4}-\eqref{4.7} and the fact that for an orthonormal frame field of ${\mathcal{D}},~\eta(e_i)=0$, for $i=1,\cdots2t$, we find
\begin{align*}
\|h\|^2&\geq4\csc^2\theta\sum_{i=1}^{2t}\sum_{r, j=1}^{s}\left(\varphi e_i(\ln f)\right)^2g(e_j^\star, e_r^\star)^2\notag\\
&+4\cot^2\theta\sum_{i=1}^{2t}\sum_{r, j=1}^{s}\left(e_i(\ln f)\right)^2g(e_j^\star, e_r^\star)^2.
\end{align*}
Hence, to satisfy \eqref{4.2}, we add and subtract the same term in the above relation and then we get
\begin{align*}
\|h\|^2&\geq4(\csc^2\theta+\cot^2%
\theta)\sum_{r, j=1}^{s}\sum_{i=1}^{2t+1}(e_i(\ln f))^2g(e_j^\star, e_r^\star)^2\notag\\
&-4(\csc^2\theta+\cot^2\theta)\sum_{r, j=1}^{s}(\xi(\ln f))^2g(e_j^\star, e_r^\star)^2.
\end{align*}
Then from \eqref{4.2} and \eqref{4.3} (i), we obtain
\begin{align*}
\|h\|^2\geq4s\left(\csc^2\theta+\cot^2\theta\right)\left(\|\nabla^T(\ln f)\|^2-1\right)
\end{align*}
which is the inequality (i). If the equality holds in (i), then from \eqref{4.11}
and \eqref{4.12}, we find
\begin{align}
\label{4.13}
h(\mathcal{D},\mathcal{D})=0,~~~h(\mathcal{D}^\theta, \mathcal{D}^\theta)=0~~{\mbox{and}}~~h(\mathcal{D}, \mathcal{D}%
_\theta)\in F\mathcal{D}^\theta.
\end{align}
If $h^\theta$ is the second fundamental form of $M_\theta$ in $M$, then we have
\begin{align}
\label{4.14}
g(h^\theta(Z, W), X)=g(\nabla_ZW, X)=-X(\ln f)\,g(Z, W)
\end{align}
 for any $X\in\Gamma(\mathcal{D})$ and $Z, W\in\Gamma(\mathcal{D}^\theta)$. Since $M_T$ is a totally geodesic submanifold in $M$ \cite{Bi,C2}, using this fact with the
first condition of \eqref{4.13}, we find that $M_T$ is totally geodesic in $\tilde M$.
Also, since $M_\theta$ is totally umbilical in $M$ \cite{Bi,C2}, using this fact with \eqref{4.14} and the second condition of \eqref{4.13}, we observe that $M_\theta$ is totally
umbilical in $\tilde M$. Moreover all conditions of \eqref{4.13} with the above fact show the minimality of $M$ in $\tilde M$. This proves the theorem completely.
\end{proof}

Now, we have the following applications of our derived results.\\

\noindent1. If we assume $\theta=\frac{\pi}{2}$ in Theorem 2 and interchange $X$ by $\varphi X$ in \eqref{4.9} for any $X\in\Gamma({\mathcal{D}}\oplus\langle\xi\rangle)$, then the warped product semi-slant submanifold becomes a contact CR-warped product of the form $M=M_T\times_fM_\perp$ in a Kenmotsu manifold. Thus, Theorem 3.4 of \cite{Kh1} is a special case of Theorem 2 as follows:

\begin{corollary} (Theorem 3.4 of \cite{Kh1}) A proper contact CR-submanifold $M$ of a Kenmotsu manifold $\tilde M$ is locally a contact CR-warped product if and only if
\begin{align*}
A_{\varphi Z}X=-\varphi X(\mu)Z,~~~~X\in\Gamma({\mathcal{D}}\oplus\langle\xi\rangle),~~~ Z\in\Gamma({\mathcal{D}}^\perp)
\end{align*}
for some smooth function $\mu$ on $M$ satisfying $W(\mu)= 0$, for each $W\in\Gamma(\mathcal{D}^\perp)$.
\end{corollary}

\noindent2. Also, if we assume $\theta=\frac{\pi}{2}$ in Theorem 3, then warped product semi-slant submanifold is of the form $M=M_T\times_fM_\perp$ i.e., $M$ becomes a contact CR-warped product. Thus, Theorem 3.1 of \cite{Ar} is a special case of Theorem 3 as below:

\begin{corollary} (Theorem 3.1 of \cite{Ar}) Let $\tilde M$ be a $(2m+1)$-dimensional Kenmotsu manifold and $M=M_1\times{_f}M_2$ be an $n$-dimensional contact CR-warped product submanifold, such that $M_1$ is a $(2t+1)$-dimensional invariant submanifold tangent to $\xi$ and $M_2$ is a $s$-dimensional anti-invariant submanifold of $\tilde M$. Then
\begin{enumerate}
\item[(i)] ~\textit{The squared norm of the second fundamental form of $M$ satisfies}
\begin{align}
\label{4.15}
\|h\|^2\geq2s\left(\|\nabla^T(\ln f)\|^2-1\right)
\end{align}
where $\nabla^T(\ln f)$ is the gradient of $\ln f$.
\item[(ii)] If the equality sign in \eqref{4.15} holds identically, then $M_1$ is a
totally geodesic submanifold and $M_2$ is a totally umbilical submanifold of $\tilde M$. Moreover, $M$ is a minimal submanifold of $\tilde M$.
\end{enumerate}
\end{corollary}

\begin{acknowledgements}
The author is thankful to the Distinguished Professor Bang-Yen Chen, Michigan State University, USA, for a nice discussion during the preparation of this article which improved the quality and presentation of this paper.
\end{acknowledgements}



\end{document}